\documentclass[11pt]{article}

\usepackage{graphicx}
\usepackage{amscd}
\usepackage[square,authoryear]{natbib}
\usepackage{marsden_article}
\newcommand{\scal}[2]{\left \langle#1,#2 \right \rangle} 
\newcommand{\pbrac}[2]{ \left \{#1,#2 \right \}}         
\newcommand{\xpbrac}[2]{\pbrac{#1}{#2}_{+}}              
\newcommand{\cpart}[2]{\frac{\partial #1}{\partial #2}}  

\begin{document}

\title{The Lie-Poisson Structure of the Euler Equations of an Ideal
Fluid}
\author{Sergiy Vasylkevych
\\ Department of Mathematics, 253-37
\\California Institute of Technology
\\ Pasadena, CA 91125
\\ \footnotesize{email: sergiy@its.caltech.edu}
\and
and\\
Jerrold E. Marsden
\\Control and Dynamical Systems Department, 107-81
\\California Institute of Technology
\\Pasadena, CA 91125
\\ \footnotesize{email: marsden@cds.caltech.edu}
}
\date{January, 2003; this version: November 30, 2003}
\maketitle

\begin{abstract}
This paper provides a precise sense in which the time $t$ map
for the Euler equations of an ideal fluid in a region in $\mathbb{R}^n$
(or a smooth compact $n$-manifold with boundary) is a Poisson map relative
to the Lie-Poisson bracket associated with the group of volume preserving
diffeomorphism group. This is interesting and nontrivial because in
Eulerian representation, the time $t$ maps need not be $C ^1$ from the
Sobolev class $H^s$ to itself (where $s > (n/2) + 1$). The idea of how
this difficulty is overcome is to exploit the fact that one does have
smoothness in the Lagrangian representation and then carefully perform a
Lie-Poisson reduction procedure.
\end{abstract}

\tableofcontents

\section{Introduction}
Hamiltonian structures play a fundamental role in mathematical
physics. It's enough to recall a few examples: classical
mechanics, electrodynamics, quantum mechanics, hydrodynamics and
general relativity. However, when applying the classical methods
and technics of symplectic geometry to PDEs, one faces significant
difficulties, both analytical and conceptual.

Part of the problem is that symplectic forms that arise in many
applications are weak symplectic forms on infinite
dimensional manifolds. More importantly, often integral curves of
PDEs are not differentiable in time in the function
spaces one would normally use; in the linear case, this corresponds to
the fact that the operators involved are unbounded. Stock examples
include the Euler and Klein-Gordon equations. When dealing with such
systems one has to pay careful attention to domains of definitions as
many standard formulas become only formal relationships. Their
justification is often cumbersome and requires some ad hoc methods.
\medskip

The goal of this paper is to contribute to the development of techniques
that are useful for the treatment of nonlinear PDEs with
non-differentiable (in time) solutions and build a framework that allows
a systematic and rigorous study of such systems and is applicable to the
broad range of physical phenomena. Previous work in this vein is
\cite{ChMa1974}.

Specifically, this article is devoted to the study of the Euler
equations for an ideal fluid on the compact manifold, the example that
provides the main inspiration and motivation. The goal is to understand
in what exact sense (if any) the flow generated by Euler's equation
consists of Poisson maps. Since the classic work of
Arnold \cite{Arnold1966}, it has been known that {\em formally} the Euler
equation could be viewed as a Hamiltonian system. (Expositions of this
may be found in \cite{ArKh1998} and \cite{MaRa1999}).

 The work of \cite{EbMa1970} showed the remarkable fact that in
appropriate function spaces, the flow of the Euler equations in
Lagrangian representation (in Sobolev function spaces $H ^s$ for $s>
(n/2) + 1$) is given by a smooth vector field and hence all the
difficulties are resolved in that context. This work also shows that one
can perform a reduction (Euler-Poincar\'e reduction) to Eulerian
representation to rigorously derive that the solutions obtained this way
satisfy the Euler equations (taking into account one derivative loss due
to the reduction procedure).

From the work of \cite{EbMa1970}, the reduced flow of the Euler equations
in $H ^s$ are known to form a continuous flow in $H ^s$ (both in time
and in the initial velocity field),  and regarded as maps from
$H ^s$ to $H ^{s - 1}$, they are $C ^1$. Another remarkable property of
the solutions also follows from this same work---namely that the
individual particle trajectories are $C ^{\infty}$ in time, a fact not
so easy to see directly in Eulerian representation (see \cite{Kato2000}).

While a version of the {\it symplectic nature} of the flow of the Euler
equations follows directly from the results in \cite{EbMa1970} (taking
into account the loss of one derivative), it is not so clear that there
is a well defined Poisson sense for the results. In fact, the work of
\cite{LeMaMoRa1986} (and many subsequent papers by other authors) shows
that in the Poisson context, this derivative loss is a nontrivial issue in
defining a good sense in which one has a Poisson manifold and in which
the Euler equations then define a Hamiltonian system in the Poisson
sense. {\it The main purposes of this paper is to fill this gap by means
of a nonsmooth Lie-Poisson reduction procedure on appropriate classes of
functions.}
\medskip

This article has the following structure. In \S\ref{s:eulersol} we give
important background information on Euler equation and manifolds of
diffeomorphisms. Then, we recall the basic ideas of Poisson reduction in
\S\ref{s:poisson}. Our results are presented in next two sections.
In \S\ref{s:weakpoisson} we prove that tangent bundle of a weak
Riemannian manifold carries a Poisson structure in an appropriate sense,
provided that the manifold possesses a {\em smooth} Riemannian
connection. The later requirement is fulfilled on the groups of
diffeomorphisms according to the work of \cite{EbMa1970}. In
\S\ref{s:eulerflow} we utilize this result to show that the flow of Euler
equation is Poisson in an appropriate sense. We conclude with short
discussion of presented results in \S\ref{s:future}.

\section{Solutions of the Euler Equation}\label{s:eulersol}
In this section we present some classical results concerning the Euler
equation that motivated our study. The notation and exposition follows
\cite{EbMa1970}.

The Euler equations on compact manifold are traditionally formulated in
the following way. Let $M$ be a compact Riemannian $n$-manifold
possibly with boundary $\partial M $. Find a time dependent
vector field
$u$, (which has an associated flow denoted $\eta_t$) such that
\begin{enumerate}
\item $u_0$ is a given initial condition with $\operatorname{div} u_0=0$
\item The Euler equations hold:
\begin{equation}
\label{euler_equations}
\cpart{u_t}{t}+\nabla_{u_t}u_t = -
\operatorname{grad} p_t
\end{equation}
 for some scalar function
$p_t:M \rightarrow \mathbb{R} $ (the pressure),
\item $\operatorname{div} u_t = 0$, and
\item $u $ is parallel to $\partial M $.
\end{enumerate}

It is standard that above equation can be formally rewritten as an ODE on
the space of divergence free vector fields with a derivative loss. But it
was discovered by \cite{EbMa1970} that this is literally true with no
derivative loss in Lagrangian representation. We recall how this
proceeds. Let
$\mu$ be a volume form on the manifold
$M$. Let
$H^s(M,N)$ denote the space of mappings of Sobolev class $s$ from
an $n$-manifold $M$ to a manifold $N$. For $s>n/2 + 1$, let
$$
{\mathcal{D}^s}=\{ \eta \in H^s(M,M) \mid \eta \ \text{is bijective
and}\ \eta^{-1} \in H^s(M,M)\} \quad \text{and}
$$
$$
\mathcal{D}^s_{\mu}=\{\eta \in \mathcal{D}^s \mid \eta^* \mu = \mu \}.
$$
Then both $\mathcal{D}^s, \mathcal{D}^s_{\mu}$ are smooth
infinite dimensional manifolds and topological groups, moreover
$\mathcal{D}^s_{\mu}$ is a closed submanifold and a subgroup of
$\mathcal{D}^s$.

Let $\widetilde\tau: T\mathcal{D}^s_{\mu} \rightarrow \mathcal{D}^s_{\mu}$
and $\tau: TM \rightarrow M$ be the canonical projections
and let $e: M \rightarrow M, \quad e(m)=m$ be the identity element of
the groups $\mathcal{D}^s_{\mu}, \mathcal{D}^s$. Then
$$
T_\eta \mathcal{D}^s = \{u \in H^s(M,TM) \mid \tau \circ u =\eta
\quad \mbox{and} \quad u \| \partial M \},
$$
$$
T_e \mathcal{D}^s_{\mu} = \mathfrak{X}_{\rm div}^s(M)=\{u \in H^s(M,TM)
\mid \tau \circ u =e, \ \operatorname{div} u =0
\quad \mbox{and} \quad u \| \partial M \},
$$
where $ \mathfrak{X}_{\rm div}^s(M)$ denotes the space of $H^s$
divergence free vector fields on $M$ that are parallel to the boundary.

A given Riemannian metric on $M$ induces a right invariant weak Riemannian
metric on $\mathcal{D}^s_{\mu}$ given by
\begin{equation}\label{e:scal}
 \scal{X}{Y}_\eta= \int_M
 \scal{X(m)}{Y(m)}_{\eta(m)} \mu(m)
 \end{equation}
for $X,Y \in T_\eta \mathcal{D}^s_{\mu}$ where scalar product
under the integral sign is taken in $M$.

As was shown in \cite{EbMa1970}, $\mathcal{D}^s_{\mu}$ possesses a
smooth Riemannian connection and, as a consequence, a smooth spray,
which we will denote $S$.

\begin{proposition} \label{t:flow} \textup{(\cite{EbMa1970})}
For $s > (n/2) + 1$, the weak Riemannian metric \textup{(\ref{e:scal})}
has a
$C ^{\infty}$ spray $S: T\mathcal{D}^s_{\mu} \rightarrow
TT\mathcal{D}^s_{\mu}$. Let $F_t: T\mathcal{D}^s_{\mu} \rightarrow
T\mathcal{D}^s_{\mu}$ be the (local, $C ^{\infty}$) flow of $S$. Let $v_t
= F_t (u_0)$ (the material velocity field) and $\eta_t= \tilde\tau (v_t)$
(the particle position field). Then the solution of the Euler equation
with initial condition $u(0)=u_0$ is given by $$ u_t = v_t \circ
\eta_t^{-1}.$$
\end{proposition}

From the properties of the diffeomorphism group, one sees that this
result shows that the Euler equations (\ref{euler_equations}) are
well-posed in $H ^s$ in Eulerian representation.

\section{Motivation: The Poisson Reduction Theorem}\label{s:poisson}
First, recall the following basic and simple result about
Poisson reduction (see, for example, \cite{MaRa1999}).

Suppose that $G$ is a Lie group that acts on a Poisson manifold
$P$ and that for each $g \in G$ the action map $\Phi_g : P
\rightarrow P$ is a Poisson map. Suppose that the quotient $P/G$
is a smooth manifold and the projection $\pi: P \rightarrow P/G$
is a submersion. Then, there is a unique Poisson structure $\{
\cdot, \cdot \}$ on $P/G$ such that $\pi$ is a Poisson map. It is
given by
$$ \{f,k\}\circ \pi = \{f \circ \pi, k \circ \pi \}_P \quad \forall k,f \in \mathcal{F}(P/G), $$
where $\{\cdot,\cdot\}_P$ is a Poisson bracket in $P$ and
$\mathcal{F}(P/G)$ is a set of smooth functions on $P/G$.

If $\mathbf{X}_{H}$ is a Hamiltonian vector field for a
G-invariant Hamiltonian $H \in \mathcal{F}(P)$, then $\pi$ also
induces reduction of dynamics. There is a function $h \in
\mathcal{F}(P/G)$ such that $H=h \circ \pi$. Since $\pi$ is a
Poisson map it transforms $\mathbf{X}_{H}$ on $P$ to $\mathbf{X}_{h}$ on $P/G$, that is, $T\pi \circ \mathbf{X}_{H} = \mathbf{X}_{h} \circ \pi$. Denoting the flow of $\mathbf{X}_{H}$ by $F_t$
and the flow of $\mathbf{X}_{h}$ by $\tilde{F_t}$ we obtain
commutative diagram
$$
\begin{CD}
P @>F_t>>P\\
@VV{\pi}V @VV{\pi}V \\
P/G @>\tilde{F_t}>> P/G
\end{CD}
$$

Our strategy is to apply the above procedure to the context of
fluids. To do so, define the map
$\pi: T\mathcal{D}^s_{\mu} \rightarrow \mathfrak{X}_{\rm div}^s$ via
$$
\pi (\eta,v)=v \circ \eta^{-1},
$$
where
$\eta \in \mathcal{D}^s_{\mu};(\eta,v) \in T_\eta \mathcal{D}^s_{\mu};
\tau \circ v =\eta$.  Let $\tilde{F_t}: \mathfrak{X}_{\rm div}^s
\rightarrow \mathfrak{X}_{\rm div}^s$ be given by
$$
\tilde{F_t} (v)= \pi \circ F_t(v)
$$
for $v \in \mathfrak{X}_{\rm
div}^s. $ By Proposition \ref{t:flow}, $\tilde{F_t}$ is the flow of
Euler equation on $\mathfrak{X}_{\rm div}^s$, i.e. $u_t = \tilde{F_t}
(u_0)$ satisfies the Euler equations (\ref{euler_equations}).

It is clear from the preceding developments that $F_t$ (as a flow of a
spray) is a flow of Hamiltonian vector field on $T \mathcal{D}^s_{\mu}$.
The following commutative diagram
$$
\begin{CD}
T\mathcal{D}^s_{\mu} @>F_t>>T\mathcal{D}^s_{\mu}\\
@VV{\pi}V @VV{\pi}V \\
T_e \mathcal{D}^s_{\mu} @>\tilde{F_t}>> T_e \mathcal{D}^s_{\mu}
\end{CD}
$$
suggests that the flow of Euler equation itself, which is obtained from
$F_t$ via Poisson reduction, should be a Hamiltonian flow in the sense
of Poisson manifolds and this is certainly formally true (see, for
instance \cite{LeMaMoRa1986} for both the case considered here as well
as the case of free boundary problems).

However, as noted in this reference and elsewhere, there are difficulties
in finding the right class of functions so that one gets a Poisson
structure in a precise sense. To justify the formal insight
in precise function spaces, one has to overcome two hurdles.

The first hurdle is that $T\mathcal{D}^s_{\mu}$ is only a weak
symplectic manifold, and therefore does not necessary carry a Poisson
bracket in any obvious way without special ad hoc hypotheses such as
``the needed functional derivatives exist'' which have long been
recognized as awkward at best.

The second hurdle is that $T\mathcal{D}^s_{\mu}$ is not a Lie group in the
usual sense (left multiplication is not smooth), and
$\pi$ is not a smooth map (inversion in $\mathcal{D}^s_{\mu}$ is
not smooth). Therefore, the well developed theory of Poisson and
Lie-Poisson reduction is not directly applicable in this case, even
though the loss of derivatives one suffers from these transformations is
well understood.

The main point of this paper is to resolve these difficulties in what
we believe is a satisfactory way. We do this in the following sections.

\section{Poisson Structures on Weak Riemannian Manifolds}
\label{s:weakpoisson}
Let $Q$ be a weak Riemannian manifold modelled on Banach space
$\mathbf{E}$ with metric $\scal{\cdot}{\cdot}$. Then ${TQ}$ possesses
a canonical {\em weak} symplectic form that is given in charts by the
following standard formula (see, e.g., \cite{MaRa1999}):
$$
\Omega(\eta,e)((e_1,e_2),(e_3,e_4))=\scal{e_1}{e_4}_{\eta}-\scal{e_2}{e_3}_{\eta}+D_\eta
\scal{e}{e_1}_{\eta} \cdot e_3 - D_\eta \scal{e}{e_3}_{\eta} \cdot
e_1,
$$
where $\eta \in Q$, $e, e_1, e_2, e_3, e_4 \in {E}$.

For a smooth function $f:M \rightarrow \mathbb{R}$ on a (strong)
symplectic manifold $(M,\Omega_1)$, let $\mathbf{X}_{f}$ denote
its Hamiltonian vector field. Then
\begin{equation}\label{f:symppoisson}
\pbrac{f}{g}=\Omega_1(X_f,X_g)
\end{equation}
makes $(M,\pbrac{\cdot}{\cdot})$ into a Poisson manifold.

Since $\Omega$ is weak, formula \ref{f:symppoisson} does not
automatically define Poisson bracket $\pbrac{f}{g}$ for arbitrary
functions $f,g \in \mathcal{F}({TQ})$ since $\mathbf{X}_{f},\mathbf{X}_{g}$ may fail to exist and even if they do, one has
to make additional hypotheses to obtain the Jacobi identity.

However, under the two additional hypothesis:
\begin{enumerate}
\item $Q$ has smooth Riemannian connection;
\item The inclusion $T_{\eta} Q \rightarrow T_{\eta}^{*} Q $ (the
literal dual space) via
$$
v(u)=\scal{v}{u}_\eta \quad \forall u \in T_\eta Q
$$
is dense,
\end{enumerate}
it will be shown that one can define a Poisson bracket on the subalgebra
$$
\mathcal{K}({TQ}) =
\left\{ f \in {TQ} \left| \cpart{f}{\eta},\cpart{f}{v} \in
C^{\infty}({TQ},{TQ}) \right\} \right.
$$
of
$\mathcal{F}({TQ})$. Here $\cpart{f}{\eta},\cpart{f}{v}$ are
{\bfi covariant partial derivatives} on ${TQ}$, the definition of which
will be given below.

This newly defined bracket {\it makes $\mathcal{K}({TQ})$ into a Lie
algebra and retains essential dynamical properties of a ``true''
Poisson bracket, including the Jacobi identity and the fact that flows of
Hamiltonian vector fields are Poisson maps and, of course, energy is
conserved.} Moreover, we will show that the bracket indeed is related to
the canonical weak symplectic form in the way that one would expect. In
the following we assume that conditions (1) and (2) are satisfied.

\paragraph{Covariant Partial Derivatives.} First, we introduce covariant
partial derivatives on
${TQ}$. Let
$\tau: {TQ} \rightarrow Q$ and $\tau_{1}: T{TQ} \rightarrow {TQ}$
be natural projections, $\Gamma: Q \supset U \times {E} \times {E}
\rightarrow {E}$ be a Christoffel map and $K: T{TQ} \rightarrow TQ$
be a connector map. In local representation,
$$
K(\eta,v,u,w)=(\eta,w+\Gamma(\eta)(v,u)).
$$
Define $\Theta:
T{TQ} \rightarrow {TQ} \bigoplus {TQ} \bigoplus {TQ}$ by
$$
\Theta = (\tau_1, T\tau,K).
$$
It is standard that $\Theta$ is a diffeomorphism (see \cite{eli67}). For
$H: {TQ} \rightarrow \mathbb{R}$ we set
$$
\cpart{H}{\eta}(V)\cdot W = dH \cdot \Theta^{-1}(V,W,0) \quad \forall \,
V,W \in T_qQ,
$$
$$
\cpart{H}{v}(V)\cdot W = dH \cdot \Theta^{-1}(V,0,W) \quad \forall \, V,W
\in T_qQ.
$$
In local representation, this reads
\begin{align*}
\cpart{H}{\eta}(\eta,v)\cdot (\eta,u) & = dH \cdot
\Theta^{-1}((\eta,v),(\eta,u),(\eta,0))
\\
& =dH\cdot(\eta,v,u,-\Gamma(\eta)(v,u)),
\end{align*}
and
\begin{align*}
\cpart{H}{v}(\eta,v)\cdot (\eta,w) & = dH \cdot
\Theta^{-1}((\eta,v),(\eta,0),(\eta,w)) \\
& =dH\cdot(\eta,v,0,w)).
\end{align*}
Similarly, for $\phi: {TQ} \rightarrow {TQ}_1$ we define
$\cpart{\phi}{\eta}, \cpart{\phi}{v}: {TQ} \rightarrow L({TQ},
T{TQ}_1)$ (here $ L({TQ}, T{TQ}_1)$ is the space of linear maps
${TQ} \rightarrow T{TQ}_1$) by
$$
\cpart{\phi}{\eta}(V)\cdot W = T\phi \cdot \Theta^{-1}(V,W,0) \quad
\forall \, V,W \in T_qQ,
$$
$$
\cpart{\phi}{v}(V)\cdot W = T\phi \cdot \Theta^{-1}(V,0,W) \quad \forall
\, V,W \in T_qQ.
$$
\bigskip

The following Lemmas are readily verified.

\begin{lemma}\label{t:prop1}
Let $X$ be a vector field on ${TQ}$, $Y$ be a vector field on
${TQ}_1$, $\phi:{TQ}_1 \rightarrow {TQ}$.  Then
\begin{align*}
dH \cdot X & = \cpart{H}{\eta} \cdot T\tau (X) + \cpart{H}{v} \cdot K
(X),\\
 \cpart{(H\circ \phi)}{\eta} \cdot Y
& = dH \cdot \left(\cpart{\phi}{\eta} \cdot Y\right), \\
  \cpart{(H\circ \phi)}{v} \cdot Y
& =dH \cdot \left( \cpart{\phi}{v} \cdot Y \right).
\end{align*}
\end{lemma}

\begin{lemma}\label{t:prop2}
For $H \in C^1(TQ,\mathbb{R})$, we have
\begin{align*}
 \cpart{H}{\eta}(\eta,v) \cdot (\eta,u) & =
\left. \frac{d}{dt}\right| _{t=0} H(\eta_t,v_t), \\
 \cpart{H}{v}(\eta,v) \cdot (\eta,w) & =
\left. \frac{d}{dt}\right|_{{t=0}}
H(\eta,v+tw),
\end{align*}
 where $(\eta_t,v_t)$ is the parallel translation of
$(\eta,v)$ along the curve $\eta_t$ with $\eta_t^{\prime} (0)=u$.
\end{lemma}

Let
\[
\mathcal{K}^k({TQ})=
\left\{ f \in C^{k+1}({TQ},\mathbb{R}) \, \left| \,
\cpart{f}{\eta},\cpart{f}{v} \in C^k({TQ},{TQ})
\right\} \right. .
\]
Now we can define the bracket
$
\{\cdot,\cdot\}$ via
\begin{equation} \label{bracket}
\{f,g\}(\eta,v)=\scal{\frac{\partial f}{\partial
\eta}(\eta,v)}{\frac{\partial g}{\partial v}{v}(\eta,v)}_\eta -
\scal{\frac{\partial f}{\partial v}(\eta,v)}{\frac{\partial
g}{\partial \eta}(\eta,v)}_\eta.
\end{equation}

\paragraph{Preliminaries on the Poisson Structure.} The following is the
first main result.

\begin{theorem}\label{t:prop4}
The bracket \textup{(\ref{bracket})} maps $\mathcal{K}^k \times
\mathcal{K}^m$ into $\mathcal{K}^{{\rm min}(k,m)-1}$ and also maps
$\mathcal{K} \times
\mathcal{K}$ into $\mathcal{K}$.
\end{theorem}

\paragraph{Remark.}
By definition of the covariant partial derivatives,  $\cpart{h}{\eta},
\cpart{h}{v}: T_\eta Q \rightarrow T^*_\eta Q$ for $h: {TQ}
\rightarrow \mathbb{R}$. The theorem asserts that if $h=\{f,g\}$
then, in fact, $\cpart{h}{\eta}(\eta,v), \cpart{h}{v}(\eta,v) \in
{TQ}$, i.e. there are $Z(\eta,v), Y(\eta,v) \in T_\eta Q$ such
that
$$ \cpart{h}{\eta} (\eta,v) \cdot X = \scal{Z}{X}, \quad  \cpart{h}{v} (\eta,v) \cdot X =
\scal{Z}{X} \quad \forall X \in T_\eta Q$$ and the maps $(\eta,v)
\rightarrow Z(\eta,v), Y(\eta,v)$ have appropriate smoothness.
\medskip

\proof Define operator $\frac{\mathbf D}{dt} = K \circ
\frac{d}{dt}$. This definition extends the usual notion of
covariant derivative from vector fields along curves on $Q$ to
arbitrary curves on ${TQ}$.
Let $f,g: {TQ} \rightarrow \mathbb{R}$ and $h=\{f,g\}$. Choosing
$(\eta_t,v_t)$ as in Lemma
\ref{t:prop2}, we obtain
\begin{align*}
& \cpart{h}{\eta}(\eta,v) \cdot (\eta,u) \\
& = \left. \frac{d}{dt}\right| _{t=0}
\scal{
\frac{\partial f}{\partial \eta} ({\eta_t},v_t)
}
{
\frac{\partial g}{\partial v} ({\eta_t},v_t)}_{\eta_t} -
\left. \frac{d}{dt} \right| _{t=0}
\scal{\frac{\partial f}{\partial v}({\eta_t},v_t)}
{\frac{\partial g}{\partial \eta}({\eta_t},v_t)}_{{\eta_t}}
\\
&
= \scal{{\left.\frac{\mathbf D}{dt}\right|_{t=0}} \frac{\partial
f}{\partial
\eta} (\eta_t,v_t)}{\frac{\partial g}{\partial v}(\eta,v)}_\eta +
 \scal{\frac{\partial f}{\partial \eta} (\eta,v)}
{ \left. \frac{\mathbf D}{dt} \right|_{t=0} \frac{\partial g}{\partial
v}({\eta_t},v_t)}_\eta
\\
& \qquad
-
\scal{
{\left. \frac{\mathbf D}{dt}\right|_{t=0}} \frac{\partial f}{\partial v}
(\eta_t,v_t)}{\frac{\partial g}{\partial \eta}(\eta,v)}_\eta -
\scal{\frac{\partial f}{\partial v}(\eta,v)}
{{\left. \frac{\mathbf D}{dt}\right|_{t=0}}
\frac{\partial g}{\partial \eta}({\eta_t},v_t)}_\eta
\\
& =
 \scal{K \cpart{}{\eta} \frac{\partial f}{\partial \eta}
(\eta,v) \cdot(\eta,u)}{\frac{\partial g}{\partial v}(\eta,v)}_\eta +
 \scal{K \cpart{}{\eta} \frac{\partial g}{\partial v}
(\eta,v) \cdot(\eta,u)}{\frac{\partial f}{\partial \eta}(\eta,v)}_\eta
\\
& \qquad
 -\scal{K \cpart{}{\eta}
\frac{\partial g}{\partial \eta} (\eta,v) \cdot(\eta,u)}
{\frac{\partial f}{\partial v}(\eta,v)} _\eta-
\scal{K \cpart{}{\eta}
\frac{\partial f}{\partial v} (\eta,v) \cdot(\eta,u)}
{\frac{\partial g}{\partial \eta}(\eta,v)}_\eta.
\end{align*}

To proceed further, we need to calculate the quantity
\[
\scal{K
\cpart{}{\eta}
\frac{\partial f}{\partial \eta} (\eta,v)
\cdot(\eta,u)}{(\eta,w)}_\eta,
\]
 where $(\eta,w)$ is an arbitrary
element of $T_\eta Q$. Let $({\eta_{ts}},{v_{ts}})$ be a
parametric surface in ${TQ}$ with the following properties:
\begin{enumerate}
\item $\left. \frac{d}{dt}\right|_{{t=0}} \eta_{t0}=u,
(\eta_{00},v_{00})=(\eta,v)$;
\item $(\eta_{t0},v_{t0})$ is a parallel translation of
$(\eta,v)$;
\item $(\eta_{t0},w_t)$ is a parallel translation of
$(\eta_{00},w_0)=(\eta,w)$;
\item $\left. \frac{d}{ds}\right|_{{s=0}} {\eta_{ts}} = w_t$ for all $s$;
\item $({\eta_{ts}},{v_{ts}})$ is a parallel translation of
$(\eta_{t0},v_{t0})$ for all $s$.
\end{enumerate}

Then, keeping in mind Lemmas \ref{t:prop1}, \ref{t:prop2} and
symmetry of Riemannian connection, one checks the following:
\begin{align*}
& \scal{K \cpart{}{\eta} \frac{\partial f}{\partial \eta} (\eta,v)
\cdot(\eta,u)}{(\eta,w)}_\eta=
 \scal{{\left. \frac{\mathbf D}{dt}\right|_{t=0}}
\frac{\partial f}{\partial \eta} (\eta_{t0},v_{t0})}{(\eta,w)}_\eta
\\
& \! \! \! \! \! \! \! \! \! \! \! \! \! \! = \frac{d}{dt}_{{t=0}}
\scal{ \frac{\partial f}{\partial \eta}
(\eta_{t0},v_{t0})}{(\eta_{t0},w_t)}_{{\eta_t}}
\end{align*}
 $$
= \left. \frac{d}{dt}\right| _{{t=0}} \frac{d}{ds}_{{s=0}}
f({\eta_{ts}},{v_{ts}})=
\frac{d}{ds}_{{s=0}} \frac{d}{dt}_{{t=0}} f({\eta_{ts}},{v_{ts}})
  = \frac{d}{ds}_{{s=0}} df \cdot
\left. \frac{d}{dt}\right|_{{t=0}}({\eta_{ts}},{v_{ts}})
 $$
 $$
= \left. \frac{d}{ds}\right|_{{s=0}} \left [\frac{\partial f}{\partial
v}(\eta_{0s},v_{0s})\cdot K \frac{d}{dt}_{{t=0}} ({\eta_{ts}},{v_{ts}})
+ \frac{\partial f}{\partial \eta} (\eta_{0s},v_{0s}) \cdot T\tau
\frac{d}{dt}_{{t=0}} ({\eta_{ts}},{v_{ts}})
 \right ]
 $$
 $$
= \left.\frac{d}{ds}\right|_{{s=0}} \left [\scal{\frac{\partial
f}{\partial v}(\eta_{0s},v_{0s})}{{\frac{\mathbf D}{dt}_{t=0}}
({\eta_{ts}},{v_{ts}})}_{\eta{0s}} +
 \scal{\frac{\partial f}{\partial \eta}
(\eta_{0s},v_{0s})}{\frac{d}{dt}_{{t=0}}
 {\eta_{ts}}}_{\eta_{0s}} \right ]
 $$
 $$
\quad  =
\scal{\frac{\mathbf D}{ds}_{s=0} \frac{\partial f}{\partial
v}(\eta_{0s},v_{0s})}{{\frac{\mathbf D}{dt}_{t=0}}
(\eta_{t0},v_{t0})}_\eta+
 \scal{\frac{\partial f}{\partial v}(\eta,v)}{\frac{\mathbf D}{ds}_{s=0}
{\left.\frac{\mathbf D}{dt}\right|_{t=0}} ({\eta_{ts}},{v_{ts}})}_\eta
 $$
 $$
 +
\scal{\frac{\mathbf D}{ds}_{s=0} \frac{\partial f}{\partial \eta}
(\eta_{0s},v_{0s})}{\frac{d}{dt}_{{t=0}} \eta_{t0}}_{\eta}+
 \scal{\frac{\partial f}{\partial \eta} (\eta,v)}
{ \left. \frac{\mathbf D}{ds}\right|_{s=0} \frac{d}{dt}_{{t=0}}
{\eta_{ts}}}_{\eta}.
 $$

\begin{lemma}\label{t:lemma1,2}
\textup{(see \cite{doc92}}). Let $\mathcal{R}$ denote the Ricci
curvature tensor. Then
$$
\frac{\mathbf D}{ds} \frac{\mathbf D}{dt} ({\eta_{ts}}, {v_{ts}})
= \frac{\mathbf D}{dt} \frac{\mathbf D}{ds}
({\eta_{ts}},{v_{ts}})+\mathcal{R}(\frac{d}{dt} {\eta_{ts}},
\frac{d}{ds} {\eta_{ts}}) ({\eta_{ts}},{v_{ts}}),
$$
$$
\frac{\mathbf D}{ds} \frac{d}{dt} {\eta_{ts}} = \frac{\mathbf
D}{dt} \frac{d}{ds} {\eta_{ts}}.
$$
\end{lemma}

By construction of $({\eta_{ts}},{v_{ts}})$, we have
${\left.\frac{\mathbf D}{dt}\right|_{t=0}} (\eta_{t0},v_{t0})=0$.
Applying lemma
\ref{t:lemma1,2} we obtain
$$
\frac{\mathbf D}{ds}_{s=0} {\left.\frac{\mathbf D}{dt}\right|_{t=0}}
({\eta_{ts}},{v_{ts}}) = \mathcal{R} \left(
\left.\frac{d}{dt}\right|_{{t=0}}
\eta_{t0},\left.\frac{d}{ds}\right|_{{s=0}}
\eta_{0s} \right) (\eta,v) = \mathcal{R}((\eta,u),(\eta,w))(\eta,v),
$$
$$
\frac{\mathbf D}{ds}_{s=0} \left.\frac{d}{dt}\right|_{{t=0}} {\eta_{ts}}
= {\frac{\mathbf D}{dt}_{t=0}} (\eta_{t0},w_{t}) = 0.
$$

Thus,
\begin{align*}
& \scal{K \cpart{}{\eta} \frac{\partial f}{\partial \eta} (\eta,v)
\cdot(\eta,u)}{(\eta,w)}_\eta \\
& = 0 +
\scal{\frac{\partial f}{\partial
v}(\eta,v)}{\mathcal{R}((\eta,u),(\eta,w))(\eta,v)}_\eta
\\
& \qquad
+\scal{K \cpart{}{\eta} \frac{\partial f}{\partial \eta} (\eta,v)
\cdot(\eta,w)}{(\eta,u)}_\eta+0
\\
&
=\scal{K \cpart{}{\eta} \frac{\partial f}{\partial \eta} (\eta,v)
 \cdot(\eta,w)}{(\eta,u)}_\eta-
\scal{\mathcal{R}((\eta,v),\frac{\partial f}{\partial
v}(\eta,v))(\eta,w)}{(\eta,u)}_\eta
\end{align*}
by Bianchi's identity. Similar calculations yield
$$
\scal{K \cpart{}{v} \frac{\partial f}{\partial v} (\eta,v)
\cdot(\eta,u)}{(\eta,w)}_\eta=
\scal{K \cpart{}{v} \frac{\partial f}{\partial v} (\eta,v)
\cdot(\eta,w)}{(\eta,u)}_\eta,
$$
$$
\scal{ K \cpart{}{\eta} \frac{\partial f}{\partial v} (\eta,v)
\cdot(\eta,u)}{(\eta,w)}_\eta
=
\scal{K \cpart{}{v} \frac{\partial f}{\partial \eta} (\eta,v)
\cdot(\eta,w)}{(\eta,u)}_\eta.
$$

Substituting this into the formulas for $\cpart{h}{\eta}$ and using
Bianchi's identity once again, we get
\begin{align*}
& \cpart{h}{\eta} (\eta,v) \cdot (\eta, u) =
\scal{K \cpart{}{\eta} \frac{\partial f}{\partial \eta} (\eta,v) \cdot
\frac{\partial g}{\partial v}(\eta,v) +
 K \cpart{}{v} \frac{\partial g}{\partial \eta}(\eta,v) \cdot
\frac{\partial f}{\partial \eta}(\eta,v)}{(\eta,u)}_\eta
\\
& -\scal{K \cpart{}{\eta} \frac{\partial g}
{\partial \eta} (\eta,v) \cdot \frac{\partial f}{\partial v}(\eta,v)
 + K \cpart{}{v} \frac{\partial f}{\partial \eta}(\eta,v) \cdot
\frac{\partial g}{\partial \eta}(\eta,v)}{(\eta,u)}_\eta\\
& \quad +
\scal{ \mathcal{R} (\frac{\partial f}{\partial v},
\frac{\partial g}{\partial \eta}) \cdot (\eta,v)} {(\eta,u)}_\eta.
\end{align*}

Similarly,
$$
 \cpart{h}{v} (\eta,v) \cdot (\eta, u) =
 \scal{K \cpart{}{\eta} \frac{\partial f}{\partial v} (\eta,v) \cdot \frac{\partial g}{\partial v}(\eta,v) +
 K \cpart{}{v} \frac{\partial g}{\partial v}(\eta,v) \cdot \frac{\partial f}{\partial \eta}(\eta,v)}{(\eta,u)}_\eta
$$
$$ -
 \scal{K \cpart{}{\eta} \frac{\partial g}{\partial v} (\eta,v) \cdot \frac{\partial f}{\partial v}(\eta,v)+
 K \cpart{}{v} \frac{\partial f}{\partial v}(\eta,v) \cdot \frac{\partial g}{\partial \eta}(\eta,v)}{(\eta,u)}_\eta.
$$

As $K$ is smooth, the statement of the theorem follows.   \qed

\paragraph{Hamiltonian Vector Fields.} The smoothness structure of
Hamiltonian vector fields is given as follows.
\begin{proposition}\label{t:prop5,6}
The vector field $\mathbf{X}_{H}$ is a $C^k$ Hamiltonian vector field
(with respect to canonical weak symplectic form) on ${TQ}$ of
class $C^k$ if and only if $H \in \mathcal{K}^k ({TQ})$.
Moreover,
\begin{equation}\label{e:ham}
\mathbf{X}_{H}(\eta,v) = \left (\eta,v,\cpart{H}{v}, -
\cpart{H}{\eta} - \Gamma(\eta)(v,\cpart{H}{v}) \right ).
\end{equation}
\end{proposition}
\proof In local representation, we have
\begin{equation}\label{e:omega}
 \Omega(\eta,e)((e_1,e_2),(e_3,e_4))=\scal{e_1}{e_4}_{\eta}-\scal{e_2}{e_3}_{\eta}+
 \scal{\Gamma(\eta)(e,e_3)}{e_1}_{\eta}-\scal{\Gamma(\eta)(e,e_1)}{e_3}_{\eta}.
\end{equation}
Indeed,
$$
 D_\eta \scal{e}{e_1}\cdot e_3 = \scal{\Gamma(\eta)(e_3,e_1)}{e}_{\eta}+
 \scal{\Gamma(\eta)(e_3,e)}{e_1}_{\eta} .
$$
Substituting this expression into the formula for $\Omega$ and
using the symmetry of $\Gamma$ we obtain the desired result.  \qed
\bigskip

Let $\mathbf{X}_{H}=(\eta,v,e_1,e_2)$ be a Hamiltonian vector
field, $Z=(\eta,v,u,w)\in T_{(\eta,v)}{TQ}$ be arbitrary. Then
$$
 \Omega(\mathbf{X}_{H},Z)=\scal{w+\Gamma(\eta)(v,u)}{e_1}-\scal{e_2+\Gamma(\eta)(v,e_1)}{u}.
$$

On the other hand, by lemma \ref{t:prop1}
$$
 \Omega(\mathbf{X}_{H},Z)= dH \cdot Z = \cpart{H}{\eta} \cdot T \tau Z +
 \cpart{H}{v} K Z=
 \cpart{H}{\eta} \cdot (\eta,u) + \cpart{H}{v}\cdot (\eta, w+
 \Gamma(u,v)).
$$

Setting $u=0$ and comparing the above expressions we see that
$\cpart{H}{v}(\eta,v) \cdot (\eta,w)= \scal{e_1}{w} \, \forall w
\in \mathbf{E}.$ Similarly, setting $w=0$ yields
\[
\cpart{H}{\eta}(\eta,v)
\cdot (\eta,u)= - \scal{e_2+\Gamma(\eta)(v,e_1)}{u} \, \forall u
\in \mathbf{E}.
\]
 Thus, $ H \in \mathcal{K}^k$.

Conversely, let $H \in \mathcal{K}^k$. Defining a vector field
$\mathbf{X}_{H}$ by formula \ref{e:ham} and substituting into
formula \ref{e:omega} one obtains for arbitrary vector $Z \in
T_{(\eta,v)} {TQ}$
$$
 \Omega(\mathbf{X}_{H},Z) = \scal{\cpart{H}{v}}{K
 Z}+\scal{\cpart{H}{\eta}}{T \tau Z}= dH \cdot Z. \qed
$$

\begin{proposition}\label{t:prop7}
 Let $f,g \in \mathcal{K}^k$ be arbitrary. Then
 $$
  \{f,g\}= \Omega(\mathbf{X}_{f},\mathbf{X}_{g}).
 $$
\end{proposition}
\proof By Proposition \ref{t:prop5,6}, the vector fields $\mathbf{X}_{f},\mathbf{X}_{g}$ are defined whenever $\{f,g\}$ is. Then
$$
 \Omega(\mathbf{X}_{f},\mathbf{X}_{g})=df \cdot \mathbf{X}_{g}=\frac{\partial f}{\partial \eta} \cdot T \tau
 \mathbf{X}_{g} + \frac{\partial f}{\partial v} K \mathbf{X}_{g} = \frac{\partial f}{\partial \eta} \cdot \frac{\partial g}{\partial v} - \frac{\partial g}{\partial \eta} \cdot \frac{\partial f}{\partial v} =
 \{f,g\}.  \qed
$$

\begin{theorem}\label{t:pbrack}
The bracket $\pbrac{\cdot}{\cdot}$ is antisymmetric, bilinear,
derivation on each factor and makes $\mathcal{K}$ into a
Lie-algebra.
\end{theorem}

\proof Antisymmetry, linearity and property of being derivation
follows directly from the definition of the bracket. By Theorem
\ref{t:prop4} $\pbrac{\cdot}{\cdot}$ leaves $\mathcal{K}$
invariant. Then, Jacobi identity follows from Proposition
\ref{t:prop7} in the usual way, for example as in \cite{MaRa1999}.
\qed
\bigskip

Now, ${TQ}$ has both symplectic and Poisson structures, and
therefore two generally different definitions of Hamiltonian
vector fields. We need to check that in our case these
coincide. To do so, let $\mathbf{X}_{f}^P$ temporarily denote the
Hamiltonian vector field with respect to Poisson structure
$\pbrac{\cdot}{\cdot}$ and $\mathbf{X}_{f}$ denotes the Hamiltonian
vector field with respect to canonical symplectic form corresponding to
function the $f$. Recall, that $\mathbf{X}_{f}^P$ is defined as a vector
field such that
$$ \mathbf{X}_{f}^P[h]=\pbrac{h}{f} \quad \forall h\in \mathcal{K}. $$
Thus, for all $h \in \mathcal{K}$,
\begin{align*}
\mathbf{X}_{f}^P [h] & = \cpart{h}{\eta} \cdot \frac{\partial f}{\partial
v} - \cpart{h}{v} \cdot \frac{\partial f}{\partial \eta}\\
& = dh \cdot
\mathbf{X}_{f}^P = \cpart{h}{\eta} \cdot T \tau
\mathbf{X}_{f}^P+\cpart{h}{v} \cdot K
\mathbf{X}_{f}^P
\end{align*}
and therefore, $T \tau \mathbf{X}_{f}^P =
\frac{\partial f}{\partial v}$ and  $K\mathbf{X}_{f}^P = -
\frac{\partial f}{\partial \eta}$. Comparing this with formula
\ref{e:ham}, we see that $\mathbf{X}_{f} \equiv \mathbf{X}_{f}^P$.
Finally, from the coordinate expression, it is easy to
see that $\mathbf{X}_{f}$ is a well defined $C^k$ vector field
for any $f \in \mathcal{K}^k$.

Previously we established that classes $\mathcal{K}^k$ are
preserved under bracketing. Unfortunately, for $f \in \mathcal{K}^k$ and a diffeomorphism $\psi: {TQ} \rightarrow {TQ}$ the
composition $f \circ \psi$ does not have to be in any class
$\mathcal{K}^m$. One can, however, compose with \em{symplectic}
diffeomorphisms.

\begin{proposition}\label{t:compdif}
Let $\psi$ be a symplectic $C^k$ diffeomorphism, $f \in \mathcal{K}^k$. Then $f \circ \psi \in \mathcal{K}^k$.
\end{proposition}
\proof We have
$$
 \mathbf{X}_{f\circ \psi} = \psi^*(\mathbf{X}_{f}),
$$
and so by Proposition \ref{t:prop5,6}, $f\circ \psi\in \mathcal{K}^k$.
\qed
\bigskip

\begin{proposition}\label{t:flowpoisson}
Let $F_t$ be a flow of a smooth Hamiltonian vector field on $TQ$.
Then $F_t$ is a Poisson, i.e. for all $f,g \in \mathcal{K}$
$$
 \pbrac{f\circ F_t}{g\circ F_t} = \pbrac{f}{g} \circ F_t.
$$
\end{proposition}
\proof $F_t$ is symplectic with respect to the weak Riemannian
form. Since $F_t$ preserves class $\mathcal{K}$, the statement
follows from Jacobi identity by the usual argument. \qed

\section{Geometric Properties of the Flow of the Euler
Equations}\label{s:eulerflow}
As we stated earlier, in \cite{EbMa1970} it
is shown that $\mathcal{D}^s_{\mu}$ carries a smooth
Riemannian connection, and therefore the results of the previous
section apply. Therefore, by those results, the space
$T \mathcal{D}^s_{\mu}$ carries a Poisson structure
(in the precise sense given there) which we denote
$\pbrac{\cdot}{\cdot}$. Let ${K}$,
$\widehat{K}$, $\widetilde{K}$ stand for the corresponding connector maps
on the underlying manifold $M$, on $\mathcal{D}^s$ and
$\mathcal{D}^s_{\mu}$ respectively, while ${\nabla}_{}$,
$\widehat{\nabla}_{}$, $\widetilde{\nabla}_{}$ are the
corresponding connections and $\Gamma$, $\widehat{\Gamma}$,
$\widetilde{\Gamma}$ are the corresponding Christoffel maps. In the
following $\scal{\cdot}{\cdot}$ denotes the Riemannian metric on $M$,
$\mathcal{D}^s$, $\mathcal{D}^s_{\mu}$ and an induced scalar
product on $\mathfrak{X}_{\rm div}^s=T_e \mathcal{D}^s_{\mu}$
depending on the context. The relationship between these metrics
is given by \ref{e:scal}.

 Recall the
notation from \S\ref{s:poisson}. Namely, let $F_t$ be the flow of
the spray on $T\mathcal{D}^s_{\mu}$, $\tilde{F_t}$ denote the
flow of Euler equation on $\mathfrak{X}_{\rm div}^s$ and $\pi :
T\mathcal{D}^s_{\mu} \rightarrow \mathfrak{X}_{\rm div}^s$, $\pi
(\eta,v) = v \circ \eta^{-1}$. Recall also that we have
the commutative diagram

\begin{proposition}\label{t:commdiag}
The following diagram is commutative:
$$
\begin{CD}
T\mathcal{D}^s_{\mu} @>F_t>>T\mathcal{D}^s_{\mu}\\
@VV{\pi}V @VV{\pi}V \\
\mathfrak{X}_{\rm div}^s @>\tilde{F_t}>> \mathfrak{X}_{\rm div}^s
\end{CD}
$$
\end{proposition}

Now we prepare and recall from \cite{EbMa1970} some useful Lemmas.

\begin{lemma}
Let $\xi \in \mathcal{D}^s_{\mu}$.
Define $R_{\xi}: \mathcal{D}^s_{\mu} \rightarrow \mathcal{D}^s_{\mu}$
via $R_{\xi} (\eta) =
\eta \circ \xi, \, \forall \xi$. Then
$$TR_{\xi} \circ F_t(v) = F_t \circ TR_{\xi} (v) \quad \forall v \in
T\mathcal{D}^s_{\mu}. $$
\end{lemma}
\proof Indeed, notice that
\begin{align*}
\frac{d}{dt} ({\eta_t} \circ \xi, \dot{{\eta_t}} \circ \xi)
& =
({\eta_t} \circ
\xi, \dot{{\eta_t}} \circ \xi, \dot{{\eta_t}} \circ \xi,
\ddot{{\eta_t}} \circ \xi) \\
&
= TTR_{\xi} ({\eta_t},\dot{{\eta_t}},\dot{{\eta_t}},
\ddot{{\eta_t}})=TTR_{\xi}S(F_t(v)) \\
& = S (TR_{\xi} F_t(v)) = S ({\eta_t} \circ \xi, \dot{{\eta_t}} \circ
\xi)
\end{align*}
 by right invariance of the spray. Thus, $TR_{\xi}F_t(v) =
({\eta_t}
\circ \xi, \dot{{\eta_t}} \circ \xi)$ is an integral curve of $S$.
Since $TR_{\xi} F_0(v) = TR_{\xi} (v)$, the statement of the Lemma
follows from uniqueness of integral curves. \qed
\bigskip

Recall that by definition, $\tilde{F_t} (V) = \pi \circ F_t(V)$
for all $V \in T_e \mathcal{D}^s_{\mu} = \mathfrak{X}_{\rm div}^s$.
Let $V = (\eta,v) \in T\mathcal{D}^s_{\mu}$. Then, using the
preceding Lemma, we obtain
\begin{align*}
\tilde{F_t} \circ \pi (V) & = \pi \circ F_t (\pi (V))
\\
&
 = \pi \circ F_t \circ TR_{\eta^{-1}}(V)= \pi \circ
 TR_{\eta^{-1}}\circ F_t(V).
\end{align*}
Notice, that
$\pi \circ TR_{\xi} = \pi$ for any $\xi \in \mathcal{D}^s_{\mu}$. Indeed,
\begin{align*}
\pi \circ TR_{\xi} (\eta,v) & = \pi (\eta \circ \xi, v \circ \xi) =
(e, v \circ \xi \circ (\eta \circ \xi)^{-1})
\\
&
= (e, v \circ \xi \circ \xi^{-1} \circ \eta^{-1})= (e,v \circ
\eta^{-1})=\pi (\eta,v).
\end{align*}
Thus $\pi \circ TR_{\eta^{-1}}=\pi$ and the Proposition is proved.
\qed
\bigskip

\paragraph{A Poisson Structure on the Lie Algebra.}
 Now, we construct a Poisson bracket $\xpbrac{\cdot}{\cdot}$
on $\mathfrak{X}_{\rm div}^s$ so that $\pi$ is a Poisson map. For
$f,g: \mathfrak{X}_{\rm div}^s \rightarrow \mathbb{R}$ such that
$df,dg: \mathfrak{X}_{\rm div}^s \rightarrow \mathfrak{X}_{\rm div}^r$
define

$$
\xpbrac{f}{g}(v)=\scal{dg(v)}{{\nabla}_{df(v)} v}-
\scal{df(v)}{{\nabla}_{dg(v)} v} .
$$

As in \S\ref{s:weakpoisson}, define
\[
\mathcal{K}^{k,s}=\{ f \in
C^{k+1}(\mathfrak{X}_{\rm div}^s, \mathbb{R}) \mid df
\in C^k(\mathfrak{X}_{\rm div}^s,\mathfrak{X}_{\rm div}^s) \}
\]
 and
\[
\mathcal{K}^{k,s}_{r,t}=\{ f \in C^{k}(\mathfrak{X}_{\rm div}^s,
\mathbb{R}) \mid df \in C^k(\mathfrak{X}_{\rm div}^r,\mathfrak{X}_{\rm
div}^t) \}.
\]

\begin{theorem}
Let $s>n/2+1$. Then $\xpbrac{\cdot}{\cdot}$ is a bilinear map
$\mathcal{K}^{k,s}\times \mathcal{K}^{k,s} \rightarrow \mathcal{K}^{k,s}_{s+1,s-1}$ and a derivation on each factor. Moreover, it
satisfies Jacobi identity on $\mathfrak{X}_{\rm div}^{s+1}$, that is
for all $f,g,h \in \mathcal{K}^{k,s}$, and $ v\in \mathfrak{X}_{\rm
div}^{s+1}$,
$$
O(v) :=\xpbrac{f}{\xpbrac{g}{h}}(v)+ \xpbrac{h}{\xpbrac{f}{g}}(v) +
\xpbrac{g}{\xpbrac{h}{f}}(v) = 0.
$$
\end{theorem}
\proof Let $f,g \in \mathcal{K}^{k,s}$. Recall, that for $r>n/2$,
$H^r(M,\mathbb{R})$ is an algebra. Thus, $(u,v) \rightarrow
{\nabla}_{u} v$ is a bilinear bounded map $\mathfrak{X}_{\rm div}^s
\times \mathfrak{X}_{\rm div}^s
\rightarrow \mathfrak{X}_{\rm div}^{s-1}$ (and $\mathfrak{X}_{\rm div}^s
\times \mathfrak{X}_{\rm div}^{s+1} \rightarrow \mathfrak{X}_{\rm
div}^{s}$), hence smooth. This implies that
$$z(v)=\xpbrac{f}{g}(v) \in C^k(\mathfrak{X}_{\rm div}^s,\mathbb{R}).$$

Bilinearity and derivation property of $\xpbrac{\cdot}{\cdot}$
trivially follows from properties of $d,{\nabla}_{}$ and
$\scal{\cdot}{\cdot}$.

Now we calculate $dz$. Let $v,u \in \mathfrak{X}_{\rm div}^{s+1}$.
Since $z \in C^k(\mathfrak{X}_{\rm div}^s,\mathbb{R})$, the Fr\'{e}chet
derivative of $z$ exists and coincides with its Gateaux
derivative. Thus, by bilinearity of scalar product and
${\nabla}_{}$,
\begin{align*}
 dz(v)\cdot u & = \left. \frac{d}{dt}\right|_{{t=0}} z(v+tu) \\
& =
 \scal{Ddg(v)\cdot u}{{\nabla}_{df(v)}v}+
 \scal{dg(v)}{{\nabla}_{Ddf(v)\cdot u} v}
\\
& \quad +
 \scal{dg(v)}{{\nabla}_{df(v)}u}-
 \scal{Ddf(v)\cdot u}{{\nabla}_{dg(v)}v} \\
& \quad -
 \scal{df(v)}{{\nabla}_{Ddg(v)\cdot u} v}-
 \scal{df(v)}{{\nabla}_{dg(v)}u}
\end{align*}

\begin{lemma}\label{t:lem5}
Let $X\in \mathfrak{X}_{\rm div}^s$, $s>n/2+1$, and let $Y,W$ be $H^s$
vector fields on M. Then
$$
 \scal{Y}{{\nabla}_{X}W}=-\scal{{\nabla}_{X}Y}{W}.
$$
\end{lemma}
\proof By the Sobolev theorems, $X$ is a $C^1$ vector field on $M$. By
properties of the Riemannian connection, for all $m \in M$
$$
 \scal{Y}{{\nabla}_{X}W}_m = - \scal{{\nabla}_{X}Y}{W}_m+X \scal{Y}{W}_m .
$$
Thus,
$$
 \scal{Y}{{\nabla}_{X}W} =-\scal{{\nabla}_{X}Y}{W} +\int_M X \scal{Y}{W}_m  \mu .
$$
Let $G_t$ be a flow of $X$ on $M$. Since $X$ is divergence free,
$\mu$ is $G_t$ invariant, i.e. $G_t^*(\mu)=\mu$, where $G_t^*$
denotes a pullback by $G_t$. Then
\begin{align*}
 \int_M X \scal{Y}{W}_m & = \int_M \frac{d}{dt}_{{t=0}}
\scal{Y}{W}_{G_t(m)} \mu \\
& =
 \frac{d}{dt}_{{t=0}} \int_M \scal{Y}{W}_{G_t(m)} G_t^*(\mu)
\\
& =
 \frac{d}{dt}_{{t=0}} \int_M G_t^*(\scal{Y}{W}_m
 \mu)
\\
& =
 \frac{d}{dt}_{{t=0}} \int_M \scal{Y}{W}_m \mu = 0. \qed
\end{align*}

\begin{lemma}\label{t:lem6}
Let $df \in C^k(\mathfrak{X}_{\rm div}^s,\mathfrak{X}_{\rm div}^t)$,
$s,t \geq 0$. Then for all $u,v,w \in \mathfrak{X}_{\rm div}^s$
$$
 \scal{Ddf(v)\cdot u}{w}=\scal{Ddf(v)\cdot w}{u}.
$$
\end{lemma}
\proof We compute as follows:
\begin{align*}
 \scal{Ddf(v)\cdot u}{w} & =
 \frac{d}{dt}_{{t=0}} \scal{df(v+tu)}{w} \\
& =
 \frac{d}{dt}_{{t=0}} \frac{d}{ds}_{{s=0}} f(v+tu+sw)
\\
& =
 \frac{d}{ds}_{{s=0}} \frac{d}{dt}_{{t=0}} f(v+tu+sw) \\
& =
 \scal{Ddf(v)\cdot w}{u}. \qed
\end{align*}

\begin{lemma}\label{t:hodge} \textup{(The Hodge Decomposition; see
\cite{EbMa1970})}. Let $X$ be an $H^s$ vector
field on
$M$, $s \geq 0$. There is an $H^{s+1}$ function $\theta$ and an
$H^s$ vector field $Y$ with $Y$ divergence free, such that
$$
 X = {\mathrm {grad}\,} \theta + Y
$$
Further, the projection maps
$$
P_e(X)=Y
$$
$$
 Q(X)={\mathrm {grad}\,} \theta
$$
are continuous linear maps on $H^s(M,TM)$. The decomposition is
orthogonal in $L^2$ sense, that is for all $Z \in \mathfrak{X}_{\rm div}^s$
\begin{equation}\label{e:hodge}
 \scal{Z}{X}=\scal{Z}{Y}=\scal{Z}{P_e X}
\end{equation}
\end{lemma}

\begin{lemma}\label{t:lem7}
There is a bilinear continuous map $B:\mathfrak{X}_{\rm div}^s \times
\mathfrak{X}_{\rm div}^{s+1} \rightarrow \mathfrak{X}_{\rm div}^s$
$(s>n/2)$ such that for all $Z \in \mathfrak{X}_{\rm div}^s, W \in
X^{s+1}, Y \in C(M,TM)$

$$
 \scal{Z}{{\nabla}_{Y}W}=\scal{B(Z,W)}{Y}
$$
\end{lemma}
\proof Fix coordinate system $\{x_i\}$ on $M$ and let $g_{ij}$
denote components of metric tensor, $Z^i$ denote components of
vector field $Z$ in the chosen system. Let
$g_{ij}g^{jk}=\delta^k_i$ (as usually, the summation on repeated
indexes is understood). Then
\begin{align*}
 \scal{Z}{{\nabla}_{Y}W} & =
 \int_M g_{ij}Z^i
\left(\cpart{W^j}{x_k}Y^k+\Gamma^j_{kr}Y^k W^r \right) \mu
\\
&=
 \int_M  g_{sm} g^{mk} g_{ij}Z^i\left(\cpart{W^j}{x_k}+\Gamma^j_{kr} W^r
\right)
 Y^s
 \mu
\\
& =
 \scal{V}{Y},
\end{align*}
where
$$
 V^m=g^{mk} g_{ij}Z^i\left(\cpart{W^j}{x_k}+\Gamma^j_{kr} W^r
 \right).
$$
Since $H^s$ is an algebra for $s>n/2$ it follows that $V$ is an
$H^s$ vector field. Now we set
$$
 B(Z,W)=P_e V
$$
and use \ref{e:hodge}. \qed
\bigskip

By Lemmata \ref{t:lem5}-\ref{t:lem7}, we have
\begin{align*}
 dz(v)\cdot u & =
 \scal{Ddg(v)\cdot P_e {\nabla}_{df(v)}v }{u}+
 \scal{Ddf(v)\cdot B(dg(v),v)}{u}+
 \scal{{\nabla}_{df(v)}dg(v)}{u}
\\
& \quad -
 \scal{Ddf(v)\cdot P_e {\nabla}_{dg(v)}v }{u}+
 \scal{Ddg(v)\cdot B(df(v),v)}{u}+
 \scal{{\nabla}_{dg(v)}df(v)}{u}.
\end{align*}

Thus for any $v \in \mathfrak{X}_{\rm div}^{s+1}$,
\begin{align*}
 d\xpbrac{f}{g}(v) & =P_e[{\nabla}_{dg(v)}df(v)-{\nabla}_{df(v)}dg(v)]
\\
& \quad +
 Ddf(v)\cdot
 B(dg(v),v)-Ddg(v)\cdot B(df(v),v)
\\
& \quad +
 Ddg(v) \cdot P_e {\nabla}_{df(v)}v -
 Ddf(v) \cdot P_e {\nabla}_{dg(v)}v,
\end{align*}
and  hence $d\xpbrac{f}{g} \in C^k(\mathfrak{X}_{\rm
div}^{s+1},\mathfrak{X}_{\rm div}^{s-1})$ and $\xpbrac{f}{g}\in
\mathcal{K}^{k,s}_{s+1,s-1}$.

\paragraph{Remark.}
If $f,g \in \mathcal{K}^{k,s} \cap \mathcal{K}^{k,s+1}$, then $v
\rightarrow P_e {\nabla}_{df(v)} dg(v)$, $v \rightarrow P_e
{\nabla}_{dg(v)} df(v)$ are
$C^k$ as maps $\mathfrak{X}_{\rm div}^{s+1} \rightarrow
\mathfrak{X}_{\rm div}^s$, hence
$\xpbrac{f}{g}\in \mathcal{K}^{k,s}_{s+1,s}$.
\medskip

Now we prove the Jacobi identity. To simplify notation, we set
$$
 B_f(v)=B(df(v),v), \quad {\nabla}_{f}(v)= P_e {\nabla}_{df(v)}v.
$$
Moreover, since in the following argument all functions are
evaluated at the same point $v\in \mathfrak{X}_{\rm div}^{s+1}$, we
will write $B_f, {\nabla}_{f},df$ instead of $B_f(v)$, etc. By
Lemmata \ref{t:lem6}-\ref{t:lem7}, we obtain
\begin{align*}
 O_{fgh}(v)
& =\xpbrac{f}{\xpbrac{g}{h}}(v) \\
& =\scal{d\xpbrac{g}{h}}{{\nabla}_{f}}-
 \scal{B_f}{d\xpbrac{g}{h}} \\
& =\scal{d\xpbrac{g}{h}}{{\nabla}_{f}-B_f} \\
& =
 \scal{P_e[{\nabla}_{dh}dg-{\nabla}_{dg}dh]}{{\nabla}_{f}-B_f}+
 \scal{Ddg \cdot (B_h - {\nabla}_{h})}{ {\nabla}_{f}-B_f} \\
& \qquad +
 \scal{Ddh \cdot ({\nabla}_{g} - B_g)}{{\nabla}_{f}-B_f} \\
& =
 \scal{[dh,dg]}{{\nabla}_{f}-B_f}+D_{ghf}-D_{hfg},
\end{align*}
where $D_{ghf}=\scal{Ddg \cdot (B_h - {\nabla}_{h})}{
{\nabla}_{f}-B_f}$ and $[\cdot,\cdot]$ is a Lie bracket of vector
fields on $M$. Notice that Lie bracket of divergence free vector
fields is divergence free.

For $s>n/2+2$
$$
\scal{[dh,dg]}{{\nabla}_{f}-B_f}=\scal{\left[\left[dh(v),dg(v)\right],df(v)\right]}{v}.
$$
Since terms of type $D_{fgh}$ cancel out in the Jacobi cycle
$$
 O(v)=O_{fgh}(v)+O_{hfg}(v)+O_{ghf}(v),
$$
and so the Jacobi identity for bracket $\xpbrac{\cdot}{\cdot}$ follows
from the Jacobi identity for vector fields. However, for $n/2+1<s
\leq n/2+2$ Lie bracket of $dh(v)$ and $dg(v)$ is an $\mathfrak{X}_{\rm div}^{s-1}$ vector field, hence merely continuous and
therefore $\left[\left[dh(v),dg(v)\right],df(v)\right]$ may fail
to exist. Therefore, in this case more care is needed.

Let
$$ A_{fgh}=\scal{df}{{\nabla}_{dg}{\nabla}_{dh}v},$$
$$
C_{fgh}=\scal{df}{{\nabla}_{[dg,df]}v}
$$
With this notation in mind, by Lemma \ref{t:lem5} and the Hodge
decomposition
$$
 \scal{[dh,dg]}{{\nabla}_{f}}=\scal{{\nabla}_{dh}dg-{\nabla}_{dg}dh}{{\nabla}_{df}v}=
 - A_{ghf}+A_{hgf}.
$$
Similarly, by definition of $B$
$$
 \scal{[dh,dg]}{B_f}=
 C_{fhg}.
$$
By a well known formula for Riemannian connection,
$$
 {\nabla}_{X}{\nabla}_{Y} Z - {\nabla}_{Y}{\nabla}_{X} Z = {\nabla}_{[X,Y]}Z,
$$
for all sufficiently smooth vector fields $X,Y,Z$. Thus,
$$
 A_{fgh}-A_{fhg}=\scal{df}{{\nabla}_{dg}{\nabla}_{dh}v - {\nabla}_{dh}{\nabla}_{df}v}=
 \scal{df}{{\nabla}_{[dg,dh]}v}=
 C_{fgh}.
$$

Thus,
$$
 \xpbrac{f}{\xpbrac{g}{h}}=-A_{ghf}+A_{hgf}-C_{fhg}+D_{ghf}-D_{hfg},
$$
and so
\begin{align*}
 & \xpbrac{f}{\xpbrac{g}{h}} +\xpbrac{h}{\xpbrac{f}{g}} +
\xpbrac{g}{\xpbrac{h}{f}}  \\
& =
 -A_{ghf}+A_{hgf}-C_{fhg}+D_{ghf}-D_{hfg}
\\
& \qquad \qquad
 -A_{fgh}+A_{gfh}-C_{hgf}+D_{fgh}-D_{ghf}
\\
& \qquad \qquad
 -A_{hfg}+A_{fhg}-C_{gfh}+D_{hfg}-D_{fgh}
\\
& =
(A_{gfh}-A_{ghf}-C_{gfh})  +(A_{fhg}-A_{fgh}
-C_{fhg})\\
& \qquad \qquad+(A_{hgf}-A_{hfg}-C_{hgf})=0.
 \qed
\end{align*}

\paragraph{Remark.}
If $df(v),dh(v) \in \mathfrak{X}_{\rm div}^s$, $s>n/2+1$, then by
Lemma \ref{t:lem5}
$$
\xpbrac{f}{h}(v)=\scal{[dh(v),df(v)]}{v}.
$$
This shows that bracket $\xpbrac{\cdot}{\cdot}$ is naturally
related to Lie-Poisson bracket on $(\mathfrak{X}_{\rm div}^s)^*$.
\medskip

Now we establish the relationship between Poisson bracket
$\xpbrac{\cdot}{\cdot}$ on $\mathfrak{X}_{\rm div}^s$ that we just
introduced and Poisson bracket $\pbrac{\cdot}{\cdot}$ on
$\mathcal{D}^s_{\mu}$. For $f,h: \mathfrak{X}_{\rm div}^s
\rightarrow \mathbb{R}$ define
$$f_R=f \circ \pi. $$

\begin{theorem}\label{t:eqdef}
Define the function spaces
$$
 C^k_r(\mathfrak{X}_{\rm div}^s)=
\left. \left\{  f \in C^k(\mathfrak{X}_{\rm div}^s,
 \mathbb{R}) \right| df(v) \in \mathfrak{X}_{\rm div}^r \forall v \in
\mathfrak{X}_{\rm div}^s \right\} ,
$$
and
$$
 C^k_r(T\mathcal{D}_{\mu}^{s})
= \left\{ f \in C^k(T\mathcal{D}_{\mu}^{s},
 \mathbb{R}) \left| \frac{\partial f}{\partial \eta}(v),\frac{\partial
f}{\partial v}(v) \in T\mathcal{D}_{\mu}^{r} \, \forall v \in
T\mathcal{D}_{\mu}^{s} \right\} \right. .
$$
Then $f_R \in C^k_r(T\mathcal{D}_{\mu}^{s+k})$ for $f\in C^k_r(\mathfrak{X}_{\rm div}^s)$ $ (r,s>n/2+1,k\geq 1)$ and for all $f,h\in
C^1_r(X^s), v\in \mathfrak{X}_{\rm div}^{s+1}$
$$
 \xpbrac{f}{h}(v)=\pbrac{f_R}{h_R}(v)=\pbrac{f\circ \pi}{g\circ
 \pi}(v).
$$
\end{theorem}
\proof Without loss of generality $s \geq r$. Since $\pi$ is not
even a $C^1$ function $\mathcal{D}^s_{\mu} \rightarrow \mathfrak{X}_{\rm div}^s$ it is not obvious that $\pbrac{f_R}{h_R}$ is defined.
However, differentiating $f_R$ and $h_R$ as functions $T\mathcal{D}_{\mu}^{s+k} \rightarrow T\mathcal{D}_{\mu}^{s} $ one obtains
the required result.

\begin{lemma}\label{t:dfrv}
Under the assumptions of the Theorem,
$$ \cpart{f_R}{v}(\eta,v)= TR_\eta df (\pi(\eta,v)). $$
\end{lemma}
\proof It is well known (\cite{EbMa1970}) that $\pi \in
C^k(T\mathcal{D}_{\mu}^{s+k},T\mathcal{D}^s_{\mu})$. Notice,
that for $(\eta,u) \in T\mathcal{D}_{\mu}^{s+k}$,
$$
\cpart{\pi}{v}(\eta,v) \cdot (\eta, u) = \frac{d}{dt}_{{t=0}}
\pi(\eta,v+tu)=(e,v \circ \eta^{-1},0,u\circ \eta^{-1})
$$
where time derivative is taken in $T\mathcal{D}^s_{\mu}$. By
lemma \ref{t:prop1}
$$
\cpart{f_R}{v}= df \cdot \widetilde{K} \cpart{\pi}{v}.
$$
Thus, by right invariance of the metric on $\mathcal{D}^s_{\mu}$
\begin{align*}
 \cpart{f_R}{v}(\eta,v) \cdot (\eta,u) & = df (v \circ \eta^{-1})
 \cdot \widetilde{K} (e,v \circ \eta^{-1},0, u \circ \eta^{-1})
\\
&
= df (v \circ \eta^{-1}) \cdot (u \circ
\eta^{-1}) \\
& =\scal{df(\pi(\eta,v)}{u\circ \eta^{-1}}_e \\
& =
\scal{TR_{\eta} df (\pi(\eta,v))}{(\eta,u)}_{\eta}. \qed
\end{align*}
\begin{lemma}\label{t:dfrn}
Under the assumptions of the Theorem
$$
 \cpart{f_R}{\eta} (\eta,v) \cdot (\eta,u)= - \scal{df (v\circ \eta^{-1})}{\widetilde{K} [T(v\circ \eta^{-1}) \circ (u \circ
 \eta^{-1})]}_e
$$
that is,
$$
  \cpart{f_R}{\eta}(\eta,v)=- TR_\eta B_f(\pi(\eta,v)).
$$
\end{lemma}
\proof First, we calculate $\cpart{\pi}{\eta}$. Let $(\eta,u) \in
T\mathcal{D}_{\mu}^{s+k}$, $({\eta_t},v_t)$ be a parallel
translation of $(\eta,v)$ with $\frac{d}{dt}_{{t=0}} {\eta_t} =u$.
Recall that
$$
\frac{d}{dt} {\eta_t}^{-1}=-T{\eta_t}^{-1} \circ \frac{d}{dt}
{\eta_t} \circ {\eta_t}^{-1}.
$$
Then, by Lemma \ref{t:prop2},
\begin{align*}
 \cpart{\pi}{\eta}(\eta,v)\cdot (\eta,u) & = \frac{d}{dt}_{{t=0}} \pi
({\eta_t},v_t) =\frac{d}{dt}_{{t=0}}
 v_t
 \circ {\eta_t}^{-1}
\\
&
=Tv_0 \circ \frac{d}{dt}_{{t=0}} {\eta_t}^{-1} +
\left( \frac{d}{dt}_{{t=0}} v_t \right) \circ \eta_0^{-1}
\\
&
=-Tv \circ T \eta^{-1} \circ u \circ \eta^{-1} -
\left( \frac{d}{dt}_{{t=0}} v_t \right) \circ \eta^{-1}.
\end{align*}

 Since connection on $\mathcal{D}^s_{\mu}$ is right invariant,
i.e.,
$$
\widetilde{K} \circ TTR_{\xi} = TR_{\xi} \circ \widetilde{K} \quad
\forall \xi \in \mathcal{D}^s_{\mu}
$$
we have
$$
 \widetilde{K} \left [\frac{d}{dt}_{{t=0}} v_t \circ \eta^{-1}\right ]=\left [\widetilde{K} \frac{d}{dt}_{{t=0}} v_t \right ] \circ \eta^{-1}
 = 0.
$$
By Lemma \ref{t:prop1}
$$
\cpart{f_R}{\eta}=df \cdot \widetilde{K} \cpart{\pi}{\eta}.
$$

Combining above equalities together, we get
$$
 \cpart{f_R}{\eta}(\eta,v) \cdot (\eta,u) = - df \cdot \widetilde{K} \left
 [T(v\circ \eta^{-1}) \circ (u \circ \eta^{-1}) \right ].
$$
$$
 =- \scal{df (v\circ \eta^{-1})}{\widetilde{K} [T(v\circ \eta^{-1}) \circ (u \circ \eta^{-1})]}_e.
$$

We claim that for all $X,Y,Z \in \mathfrak{X}_{\rm div}^s$
\begin{equation}\label{e:kdsmcman}
 \scal{Z}{\widetilde{K} [TX \circ Y]}=\scal{Z}{{\nabla}_{Y}X}.
\end{equation}
Recall that by construction (see \cite{EbMa1970}),
$$
 \widetilde{K}=P \circ \widehat{K},
$$
$$
 P=TR_\eta \circ P_e \circ TR_\eta^{-1},
$$
$$
 \widehat{K}(Y) = {K} \circ Y,
$$
By a well known formula of differential geometry, we have
$$
 K \circ TX \circ Y={\nabla}_{Y}X,
$$
and hence
$$
 \widetilde{K}[TX \circ Y] = P_e[{\nabla}_{Y}X].
$$
By the Hodge decomposition
$$
 \scal{Z}{\widetilde{K}[TX \circ Y]} =
 \scal{Z}{{\nabla}_{Y}X}=\scal{B(Z,X)}{Y}.
$$

By the above developments and right invariance of metric on
$\mathcal{D}^s_{\mu}$, we have
$$
\cpart{f_R}{\eta}(\eta,v) \cdot (\eta,u) = -\scal{B_f(v \circ
\eta^{-1})}{u \circ \eta^{-1}}=-\scal{TR_\eta B_f
(\pi(\eta,v))}{u}_\eta. \qed
$$

Calculating
$\pbrac{f_R}{h_R}$ at $v \in \mathfrak{X}_{\rm div}^{s+1}$ by Lemmata
\ref{t:dfrv},\ref{t:dfrn}, we obtain
\begin{align*}
 \pbrac{f_R}{h_R}(v)
& =-
 \scal{B_f(v)}{dh(v)}+\scal{B_h(v)}{df(v)} \\
& =
 - \scal{df(v)}{{\nabla}_{dh(v)}v}+\scal{dh(v)}{{\nabla}_{df(v)}v}
\\
& =
 \xpbrac{df}{dh}(v). \qed
\end{align*}

\begin{proposition}\label{t:pipoisson}
Map $\pi : T\mathcal{D}^s_{\mu} \rightarrow \mathfrak{X}_{\rm div}^s$ is a Poisson map, i.e. for all $f,h \in C^{1}_{r}(\mathfrak{X}_{\rm div}^s)$ pointwise in $T\mathcal{D}_{\mu}^{s+1}$ $
(r,s>n/2+1)$
$$ \pbrac{f\circ \pi}{h\circ \pi}=\xpbrac{f}{h}\circ \pi.$$
\end{proposition}
\proof Since $\pi$ is the identity on $\mathfrak{X}_{\rm div}^s$, the
statement follows immediately from Theorem \ref{t:eqdef}. \qed

\begin{proposition}\label{t:fpoisson}
Let $v \in T\mathcal{D}_{\mu}^{r}$ and $f,g \in C^1(T\mathcal{D}_{\mu}^{r} , \mathbb{R})$ are such that $\frac{\partial
f}{\partial v}(F_t(v))$ , $\frac{\partial f}{\partial
\eta}(F_t(v))$, $\frac{\partial g}{\partial v}(F_t(v))$,
$\frac{\partial g}{\partial \eta}(F_t(v))
\in T\mathcal{D}^s_{\mu}$, $r,s>n/2+1$. Then

$$ \pbrac{f \circ F_t}{g \circ F_t}(v) = \pbrac{f}{g} ( F_t
(v)). $$ In particular, $F_t$ preserves $C^1_s(T\mathcal{D}_{\mu}^{r})$ and for $f,h \in \mathcal{K}^{1,s}$ pointwise in
$T\mathcal{D}_{\mu}^{s+1}$
$$
 \pbrac{f \circ \pi \circ F_t}{h \circ \pi \circ F_t}(v) = \pbrac{f \circ
 \pi}{h \circ \pi} ( F_t (v)).
$$
\end{proposition}
\proof Without loss of generality $r \geq s$. First, we notice
that covariant partial derivatives of $f\circ F_t,g \circ F_t$ at
$v$ are elements of $T\mathcal{D}^s_{\mu}$. Indeed,
$$ \frac{\partial }{\partial \eta}(g \circ F_t) (v) \cdot u = \scal{\frac{\partial g}{\partial \eta}(F_t(v))}{T \widetilde \tau \cpart{F_t}{\eta}(v) \cdot u} +
 \scal{\frac{\partial g}{\partial v}(F_t(v))}{\widetilde K \cpart{F_t}{\eta} (v) \cdot u}.
$$
There is a function $\tilde{g} \in \mathcal{K}(T\mathcal{D}^s_{\mu})$ such that
$$
 \frac{\partial g}{\partial v}(F_t(v))=\cpart{\tilde{g}}{v}(F_t(v)), \quad
 \frac{\partial g}{\partial \eta}(F_t(v))=\cpart{\tilde{g}}{\eta}(F_t(v)).
$$
Thus,
$$
 \frac{\partial }{\partial \eta}(g \circ F_t) (v) \cdot u = \frac{\partial }{\partial \eta}(\tilde{g} \circ F_t) (v) \cdot
 u.
$$
However, by Proposition \ref{t:compdif} $\tilde{g} \circ
F_t \in \mathcal{K}(T\mathcal{D}^s_{\mu})$ for any $\tilde{g} \in
\mathcal{K}(T\mathcal{D}^s_{\mu})$, hence
there is $Z_g \in C^\infty(T\mathcal{D}^s_{\mu},T \mathcal{D}^s_{\mu})$
such that for all $u $,
$$
\frac{\partial }{\partial \eta} (g \circ F_t) (v) \cdot u=
\frac{\partial }{\partial \eta} (\tilde{g} \circ F_t) (v) \cdot u =
\scal{Z_g(v)}{u}.
$$
In a similar sense, one shows that
$\frac{\partial }{\partial v} (f \circ F_t)(v) \in T \mathcal{D}^s_{\mu}$.

Thus, $\pbrac{f \circ F_t }{g \circ F_t}(v)$ is well defined and
depends only on values of $\frac{\partial f}{\partial v},
\frac{\partial f}{\partial \eta},
 \frac{\partial g}{\partial \eta}, \frac{\partial g}{\partial v}$ calculated at point $F_t(v)$.
 However, $\pbrac{f}{g} \circ F_t(v)$ also depends only on values
 of covariant partial derivatives at $F_t(v)$. Then, we choose
 $\tilde{f},\tilde{g} \in \mathcal{K}(T\mathcal{D}^s_{\mu})$ such that
\begin{align*}
\frac{\partial f}{\partial
v}(F_t(v)) & =\cpart{\tilde{f}}{v}(F_t(v)) \\
\frac{\partial f}{\partial
\eta}(F_t(v)) & =\cpart{\tilde{f}}{\eta}(F_t(v)),\\
\frac{\partial g}{\partial v}(F_t(v)) & =\cpart{\tilde{g}}{v}(F_t(v)),\\
\frac{\partial g}{\partial
\eta}(F_t(v)) & =\cpart{\tilde{g}}{\eta}(F_t(v)).
\end{align*}
The equality
$$
\pbrac{ \tilde{f} \circ F_t}{\tilde{g} \circ F_t}(v) =
\pbrac{\tilde{f}}{\tilde{g}} \circ F_t(v)
$$
follows from Proposition \ref{t:flowpoisson}. By the preceding arguments,
the same holds if we replace $\tilde{f},\tilde{g}$ with
$f,g$. This concludes the first part of the Proposition. The
second part then follows. \qed

\begin{theorem}\label{s:redpoisson}
The map $ \tilde{F_t}$ is Poisson with respect to the bracket
$\xpbrac{\cdot}{\cdot}$.
\end{theorem}
\proof Let $f,h \in \mathcal{K}^{k,s}$. Then $f \circ \pi \in
C^1_s(\mathfrak{X}_{\rm div}^{s+1})$. By Proposition \ref{t:fpoisson}.
$$f \circ \tilde{F_t}= f \circ \pi \circ F_t \in C^1_s(\mathfrak{X}_{\rm div}^{s+1})$$ and we have pointwise in
$\mathfrak{X}_{\rm div}^{s+2}$:
\begin{align*}
 & \xpbrac{f \circ \tilde{F_t}}{ h \circ \tilde{F_t}} \quad
\hbox{(Theorem
\ref{t:eqdef})} \\
& \qquad  = \pbrac{ f \circ \tilde{F_t} \circ \pi}{h \circ \tilde{F_t}
\circ \pi}
\quad \hbox{(Proposition \ref{t:commdiag} )}\\
& \qquad = \pbrac{ f \circ \pi \circ F_t}{h \circ \pi \circ {F_t}} \quad
\hbox{ (Proposition \ref{t:fpoisson})} \\
& \qquad = \pbrac{ f \circ \pi}{h \circ \pi} \circ F_t \quad
\hbox{(Proposition \ref{t:pipoisson})} \\
& \qquad = \xpbrac{f}{h} \circ \pi \circ {F_t} = \xpbrac{f}{h} \circ
\tilde{F_t}. \qed
\end{align*}

\section{Conclusions}\label{s:future}
In the previous sections we successfully implemented a nonsmooth
Lie-Poisson reduction technique for the study of the Euler equations of
ideal fluid flow. This enabled us to find a precise sense in which the
flow of Euler equation on the Lie algebra of divergence free vector
fields (parallel to the boundary of the fluid region) is a Hamiltonian
system in the Poisson sense and that the flow consists of Poisson maps,
despite the fact that this flow is believed (as maps from
$H^s$ to $H^s$) to be continuous, but not differentiable.

A key part of this process was to introduce a Poisson structure on the
space of divergence free vector fields. As one
would expect from the bracket derived via a type of Lie-Poisson
reduction, this bracket is closely related to the formal Lie-Poisson
bracket on the dual to the Lie algebra of divergence free vector fields.

Even though we consider only Euler's equation, the technique
developed here is directly applicable to several other important
systems---those which can be written as an ODE on groups of
diffeomorphisms, such as the following:

\begin{enumerate}
\item The Camassa-Holm (CH) equation on $S^1$---see \cite{CaHo1993}:
$$
u_t - u_{txx} = -3uu_x + 2u_xu_{xx} + uu_{xxx}.
$$
\item The averaged Euler equations (or the LAE-$\alpha$ equations)---see
\cite{HoMaRa1998a, HoMaRa1998b}:
$$
\partial_t (1-\alpha^2\Delta) u + (u
\cdot\nabla)  (1-\alpha^2\Delta)u
-\alpha^2 (\nabla u )^T
\cdot
\Delta  u = - \operatorname{grad}p\,,
$$
where $\operatorname{div} u = 0$ and $u$ satisfies appropriate
boundary conditions, such as the no-slip conditions $u = 0$ on $\partial
M $.
\item The EPDiff equation (also called the averaged template matching
equation) on a compact manifold $M$---see \cite{HoMa2003} and
\cite{HiMaAr2001}:
\begin{align*}
 & u_t - \alpha^2 \Delta u + u (\operatorname{div} u) -
\alpha^2 (\operatorname{div} u) \Delta u
  + (u \cdot \nabla) u  \\
& \qquad \qquad - \alpha^2(u\cdot \nabla) \Delta u +
(Du)^T\cdot u - \alpha^2(Du)^T \cdot \Delta u = 0,
\end{align*}
with appropriate boundary conditions, such as the no-slip
conditions $u = 0$ on $\partial M $. The EPDiff equations
reduce to the CH equations in the case $M = S^1$.
\end{enumerate}

These equations may be derived as the right reduction to the identity
of the geodesic motion on the appropriate Lie group
(see, for example, \cite{CaHo1993} and \cite{mis98,mis02} for the case
of the CH equations), and the preceding references for the other
equations. The crucial technical fact that enables our methods to
work in both cases is the smoothness of the spray on the Lie group. For
the case of the CH equations and the LAE-$\alpha$ equations on regions
with no boundary, this is due to
\cite{Shkoller1998} and for regions with boundary to \cite{MaRaSh2000}.
For the case of the EPDiff equations, a rather convincing plausibility
argument is given \cite{HoMa2003}.

One important direction in which we would like to pursue these
ideas is that of nonsmooth solutions. Even for the ideal Euler
equations, this is interesting because of the singular solutions, such
as point vortices, vortex filaments and sheets. They clearly have
themselves an interesting Poisson structure, as was investigated by
\cite{MaWe1983} and \cite{LaPe1991}. There are similar interesting
singular solutions for the EPDiff equations, whose geometry is
investigated in \cite{HoMa2003}. It would be very interesting if, on the
smaller spaces appropriate for these classes of singular solutions that
are introduced in these references, the smooth spray property still holds
and, if that is the case, whether or not one could then carry out the
program in the present paper.

Another interesting direction for the present research is to the case of
free boundary problems, a notoriously difficult case for infinite
dimensional  Poisson structures, even at the formal level (see
\cite{LeMaMoRa1986}, \cite{KrMaSc1993}, \cite{KrMaMa1999} and
\cite{Bering2000}.)


\begin{thebibliography}{}


\bibitem[Arnold(1966)]{Arnold1966}
Arnold, V.~I. [1966],
Sur la g\'{e}om\'{e}trie differentielle des groupes de {L}ie de dimenson
infinie et ses applications \`{a} l'hydrodynamique des fluids parfaits,
{\em Ann. Inst. Fourier, Grenoble} \textbf{16}, 319--361.


\bibitem[Arnold and Khesin(1998)]{ArKh1998}
Arnold, V.~I. and B.~Khesin [1998],
{\em Topological Methods in Hydrodynamics},
  volume 125 of {\em Appl. Math. Sciences}.
\newblock Springer-Verlag.



\bibitem[Bering(2000)]{Bering2000}
Bering, K. [2000],
Putting an edge on the Poisson bracket,
{\em J. Math. Phys.} \textbf{41}, 7468--7500.



\bibitem[Camassa and Holm(1993)]{CaHo1993}
Camassa, R. and D.~D. Holm [1993],
An integrable shallow water equation with
  peaked solitons, {\em Phys. Rev. Lett.} \textbf{71}, 1661--1664.

\bibitem[Chernoff and Marsden(1974)]{ChMa1974}
Chernoff, P.~R. and J.~E. Marsden [1974], {\em Properties of Infinite
  Dimensional {H}amiltonian systems}, volume 425 of {\em Lecture Notes in
  Math.}
\newblock Springer, New York.


\bibitem[Do Carmo(1992)]{doc92}
Do Carmo, M. [1992], \emph{Rimannian geometry}, Birkhauser, Boston.


\bibitem[Ebin and Marsden(1970)]{EbMa1970}
Ebin, D.~G. and J.~E. Marsden [1970], Groups of diffeomorphisms and the motion
  of an incompressible fluid, {\em Ann. of Math.} \textbf{92}, 102--163.

\bibitem[Eliasson(1967)]{eli67}
Eliasson, H. [1967], Geometry of manifolds of maps, {\it J.
Differential Geometry}
  \textbf{1}, 169--194.


\bibitem[Hirani et~al.(2001)Hirani, Marsden, and Arvo]{HiMaAr2001}
Hirani, A., J.~E. Marsden, and J.~Arvo [2001], Averaged template matching
  equations, {\em Springer Lecture Notes in Computer Science} \textbf{2134},
  528--543.

\bibitem[Holm and Marsden(2003)]{HoMa2003}
Holm, D.D. and J.~E. Marsden [2003], Momentum maps and measure valued
solutions (peakons, filaments, and sheets) of the Euler-Poincar\'{e}
equations for the diffeomorphism group. \newblock In Marsden, J.~E. and
T.~S. Ratiu, editors, {\em Festshrift for Alan
  Weinstein (to appear)}. Birkh\"auser Boston.


\bibitem[Holm et~al.(1998a)Holm, Marsden, and Ratiu]{HoMaRa1998a}
Holm, D.~D., J.~E. Marsden, and T.~S. Ratiu [1998a],
{E}uler--{P}oincar\'e models of ideal fluids with nonlinear dispersion,
{\em Phys. Rev. Lett.} \textbf{349}, 4173--4177.


\bibitem[Holm et~al.(1998b)Holm, Marsden, and Ratiu]{HoMaRa1998b}
Holm, D.~D., J.~E. Marsden, and T.~S. Ratiu [1998b], The
{E}uler--{P}oincar\'{e}
  equations and semidirect products with applications to continuum theories,
  {\em Adv. in Math.} \textbf{137}, 1--81.



\bibitem[Kato(2000)]{Kato2000}
Kato, T. [2000], On the smoothness of trajectories in incompressible perfect
  fluids, {\em Cont. Math., Am. Math. Soc.} \textbf{263}, 109--130.


\bibitem[Kruse et~al.(1999)Kruse, Mahalov, and Marsden]{KrMaMa1999}
Kruse, H.~P., A.~Mahalov, and J.~E. Marsden [1999], On the Hamiltonian
  structure and three-dimensional instabilities of rotating liquid bridges,
  {\em Fluid Dyn. Research} \textbf{24}, 37--59.


\bibitem[Kruse et~al.(1993)Kruse, Marsden, and Scheurle]{KrMaSc1993}
Kruse, H.~P., J.~E. Marsden, and J.~Scheurle [1993], On uniformly rotating
  fluid drops trapped between two parallel plates, {\em Lect. in Appl. Math.,
  AMS} \textbf{29}, 307--317.


\bibitem[Langer and Perline(1991)]{LaPe1991}
Langer, J. and R.~Perline [1991], Poisson geometry of the filament equation,
  {\em J. of Nonlinear Sci.} \textbf{1}, 71--94.

\bibitem[Lewis et~al.(1986)Lewis, Marsden, Montgomery, and Ratiu]
{LeMaMoRa1986}
Lewis, D., J.~E. Marsden, R.~Montgomery, and T.~S. Ratiu [1986], The
  {H}amiltonian structure for dynamic free boundary problems,
{\em Physica D}  \textbf{18}, 391--404.


\bibitem[Marsden and Ratiu(1999)]{MaRa1999}
Marsden, J.~E. and T.~S. Ratiu [1999], {\em Introduction to Mechanics and
Symmetry}, volume~17 of {\em Texts in Applied Mathematics, vol. 17; 1994,
Second Edition, 1999}. \newblock Springer-Verlag.


\bibitem[Marsden et~al.(2000)Marsden, Ratiu, and Shkoller]{MaRaSh2000}
Marsden, J.~E., T.~Ratiu, and S.~Shkoller [2000], The geometry and analysis of
  the averaged Euler equations and a new diffeomorphism group, {\em Geom.
  Funct. Anal.} \textbf{10}, 582--599.


\bibitem[Marsden and Weinstein(1983)]{MaWe1983}
Marsden, J.~E. and A.~Weinstein [1983], Coadjoint orbits, vortices and
  {C}lebsch variables for incompressible fluids, {\em Physica D} \textbf{7},
  305--323.

\bibitem[Misiolek(1998)]{mis98}
Misiolek, G. [1998], A shallow water equation as a geodesic
flow on the Bott-Virasoro group, {\it J. Geom. Phys.} \textbf{24},
203--208.

\bibitem[Misiolek(2002)]{mis02}
Misiolek, G. [2002], Classical solutions of the periodic
camassa-holm equation,
  {\it Geom. Func. Anal.} \textbf{12}, 1080--1104.




\bibitem[Shkoller(1998)]{Shkoller1998}
Shkoller, S. [1998], Geometry and curvature of diffeomorphism groups with $H^1$
  metric and mean hydrodynamics, {\em J. Funct. An.} \textbf{160}, 337--365.

\end{thebibliography}
\end{document}